# CHARACTERIZATION OF INVARIANT MEASURES AT THE LEADING EDGE FOR COMPETING PARTICLE SYSTEMS

By Anastasia Ruzmaikina and Michael Aizenman[1]

*University of Virginia and Princeton University*

We study systems of particles on a line which have a maximum, are locally finite and evolve with independent increments. "Quasi-stationary states" are defined as probability measures, on the $\sigma$-algebra generated by the gap variables, for which joint distribution of gaps between particles is invariant under the time evolution. Examples are provided by Poisson processes with densities of the form $\rho(dx) = e^{-sx} s\,dx$, with $s > 0$, and linear superpositions of such measures. We show that, conversely, any quasi-stationary state for the independent dynamics, with an exponentially bounded integrated density of particles, corresponds to a superposition of Poisson processes with densities $\rho(dx) = e^{-sx} s\,dx$ with $s > 0$, restricted to the relevant $\sigma$-algebra. Among the systems for which this question is of some relevance are spin-glass models of statistical mechanics, where the point process represents the collection of the free energies of distinct "pure states," the time evolution corresponds to the addition of a spin variable and the Poisson measures described above correspond to the so-called REM states.

**1. Introduction.** Competitions involving large numbers of contestants are an object of interest in various fields. One could list here the energy levels of complex systems and the free energies of competing extremal states of spin-glass models [10] and include a broad range of other examples. We are particularly interested in dynamical situations where the competition continues in "time," though time may be interpreted loosely. For example, in the motivating example of spin-glass models [10], a point process on the line represents the collection of the free energies of distinct "pure states" of a system of many spin variables, and the "time evolution" corresponds to the incorporation in the system of yet another spin variable.

---

Received February 2003; revised February 2004.
[1]Supported in part by NSF Grant PHY-99-71149.

Influenced by the terminology of statistical mechanics, we use here the term *state* to mean a probability measure on the relevant $\sigma$-algebra of subsets of the space of the point process configurations. For much of the discussion which follows, the relevance would be limited to the information concerning only the relative positions of the points, relative to the one which leads at the given instant.

As in the pictures seen in marathon races, often the point process describing the relative positions appears to be time invariant. We refer to such states as quasi-stationary.

In this paper we characterize the quasi-stationary states for the class of systems in which the evolution occurs by independent identically distributed increments of the individual contestants. The main result is that any such state, of a point process with locally finite configurations with more than one point and exponentially bounded density, corresponds to a linear superposition of Poisson processes with densities of the form

$$\rho(dx) = e^{-sx} s\, dx \tag{1.1}$$

with $s > 0$. This may be rephrased by saying, in the terminology coined by Ruelle [11] (who invokes the work of Derrida [6]), that all quasi-stationary states correspond to superpositions of the random energy model (REM) states.

Remark. Our main result may have a familiar ring to it, since the above distributions are known to describe the "Type-I" case of the extremal statistics [8].

Remark. It would be of interest to see an extension of the classification of the quasi-stationarity to a broader class of dynamics where the evolution may exhibit correlations. One may note that the REM states have an extension, based on a hierarchical construction, to the family of the GREM states [2, 11], which exhibit quasi-stationarity under a broad class of correlated dynamics. Is that structure singled out in some way by its broader quasi-stationary properties?

In the following section we introduce the concepts more explicitly. We refer to the system as the Indy-500 model, ignoring the fact that for a number of obvious reasons this is not a faithful description of the dynamics in that well-known car race.

## 2. The Indy-500 model.

The configuration space of the Indy-500 model is the space $\Omega$ of infinite configurations of points on the line, which are locally finite and have a maximum (in the order of $\mathbb{R}$). Its elements $\omega \in \Omega$ can also be described as sequences, $\omega = \{x_n\}_{n=1,2,\dots}$, with

$$x_1 \geq x_2 \geq \cdots \quad \text{and} \quad x_n \to -\infty. \tag{2.1}$$

(Variables written as $x_n$ should by default be understood to be ordered.) In the time evolution considered here the points evolve by independent increments.

As is generally the case with stochastic evolutions, the dynamics can be presented in two ways: as a stochastic map, in which the configuration $\omega \in \Omega$ is changed in a random way—through the independent increments, or as a reversible transformation taking place in a larger space, which encompasses the full information about both the future and the past dynamics. Our terminology is based on the former view; however, the second perspective provides a useful basis for the intuition guiding the analysis.

Thus, the time evolution is given by a stochastic map determined by the collection of random variables $\eta = \{h_n\}_{n=1,2,\dots}$:

$$T_\eta : \{x_n\} \mapsto \{\tilde{x}_n\} \qquad \text{with } \tilde{x}_n = x_{\Pi_n} + h_{\Pi_n}, \tag{2.2}$$

where $h_n$ are independent random variables with a common probability distribution $g(dh)$ on $\mathbb{R}$, and $\Pi$ is a permutation of $\mathbb{N}$, which depends on both $\omega$ and $\eta$, aimed at recovering the monotonicity for $\tilde{x}_n$. In other words, $\Pi = \Pi(\omega, \eta)$ is a relabeling of the moving particles according to the new order.

For a given probability measure $\mu(d\omega)$ on $\Omega$, we denote by $T\mu$ the corresponding probability distribution of the one-step evolved configuration $\{\tilde{x}_\alpha\}$. To be explicit: the average over $T\mu$ corresponds to averaging over both $\mu$ and $\eta$.

One needs to pay some attention to the $\sigma$-algebras on which the measures $\mu$ and $T\mu$ are to be defined. Since we are interested in the classification of states which are only *quasi-stationary*, we allow those to correspond to probability measures defined on a smaller $\sigma$-algebra than the one usually used for point processes on a line. (Such a change makes the result only stronger.)

The standard $\sigma$-algebra, which is natural for the state space of particle configurations, is generated by the occupation numbers of finite intervals (see, e.g., [4]). Let us denote it by $\mathcal{B}$. Measurable functions include all $\psi : \Omega \to \mathbb{R}$ of the form

$$\psi_f(\omega) := \sum_n f(x_n) \tag{2.3}$$

with bounded measurable functions $f : \mathbb{R} \to \mathbb{R}$, of compact support. However, in this work we are interested in probability measures on the smaller $\sigma$-algebra $\widetilde{\mathcal{B}}$ generated by functions which are invariant, $\psi(S_b\omega) = \psi(\omega)$, under the uniform shifts

$$S_b : \{x_n\} \mapsto \{\tilde{x}_n\} \qquad \text{with } \tilde{x}_n = x_n + b. \tag{2.4}$$

Functions which are measurable with respect to $\widetilde{\mathcal{B}}$ depend on the configuration only through the sequence of the distances of the particles from the leading one:

$$u_n = x_1 - x_n. \tag{2.5}$$

Thus, a probability measure $\mu$ on $(\Omega, \widetilde{\mathcal{B}})$ is uniquely determined by the "modified probability generating functional" (MPGFL)

$$\widetilde{G}_\mu(f) = \mathbb{E}_\mu\left(\exp\left\{-\sum_n f(x_1 - x_n)\right\}\right), \tag{2.6}$$

with $f(\cdot)$ ranging over smooth positive functions of compact support. [The regular "probability generating functional" is defined without $x_1$ in (2.6).]

One can now formulate a number of distinct "steady-state" conditions, where the term *state* refers to a probability measure on a suitable $\sigma$-algebra, which is not always the same.

Definition. A *stationary state* is a probability measure $\mu(d\omega)$ on $(\Omega, \mathcal{B})$ which is invariant under the stochastic map $T$, that is, $T\mu = \mu$, or more explicitly,

$$\mathbb{E}_\mu(\psi(T\omega)) = \mathbb{E}_\mu(\psi(\omega)) \tag{2.7}$$

for any $\mathcal{B}$-measurable $\psi$, where the expectation functional $\mathbb{E}_\mu$ includes an average over both $\omega$ (distributed by $\mu$) and $T$ [determined through $\{h_n\}$, as in (2.2)].

A *steady state* is a probability measure $\mu(d\omega)$ on $(\Omega, \mathcal{B})$ for which there is a nonrandom $V$ (= the "front velocity") such that $T\mu = S_V\mu$, that is,

$$\mathbb{E}_\mu(\psi(T\{x_n\})) = \mathbb{E}_\mu(\psi(\{x_n + V\})) \tag{2.8}$$

for all $\mathcal{B}$-measurable functions $\psi$.

A *quasi-stationary* state is a probability measure $\mu(d\omega)$ on the $\sigma$-algebra $\widetilde{\mathcal{B}}$ (sub-$\sigma$-algebra of $\mathcal{B}$) such that (2.7) restricted to *shift-invariant* functions $\psi$ holds, that is, for which

$$\mathbb{E}_\mu(\psi(\{u_n\})) = \mathbb{E}_\mu(\psi(\{\tilde{u}_n\})) \tag{2.9}$$

with $\{u_n\}$ the gaps defined by (2.5), and $\{\tilde{u}_n\}$ the gaps for the configuration $\tilde{\omega} = T\omega$.

For an alternative characterization of quasi-stationary measures, in terms which are more standard for point processes, let us note that each configuration is shift-equivalent to a unique element of the set

$$\Omega_o = \{\{x_n\} | x_1 = 0\}. \tag{2.10}$$

The "normalizing shift" $S : \omega \mapsto S_{-x_1(\omega)}\omega$ induces a measurable map from $(\Omega, \widetilde{\mathcal{B}})$ to $(\Omega_o, \mathcal{B}) \subset (\Omega, \mathcal{B})$, and thus also a map (for which we keep the symbol $S$) which associates to each probability measure $\mu$ on $(\Omega, \widetilde{\mathcal{B}})$ a probability measure $S\mu$ on $(\Omega, \mathcal{B})$, supported on $\Omega_o$. The measure $\mu$ is quasi-stationary if and only if the corresponding measure $S\mu$ is invariant under $ST$—the time evolution followed by the normalizing shift.

Stationarity is a special case of the steady state, and the latter reduces to it when viewed from a frame moving at a fixed speed. Quasi-stationarity is the less demanding property of the three mentioned above, and is the condition of interest if one follows only the relative positions.

Through a combination of the results in [9] and [3] one may conclude that any *steady* state of the Indy-500 model, whose jump distribution satisfies the *nonlattice* condition (meaning that its support is not contained in any set of the form $a + b\mathbb{Z} \subset \mathbb{R}$), is a Poisson process with a density of the form $\rho(dx) = se^{-sx}\,dx$. These are the REM states which are discussed in the next section.

Our main result is that for the infinite systems discussed here *quasi-stationary* probability measures can only be linear superpositions (as probability measures) of the above steady states restricted to $\widetilde{\mathcal{B}}$.

REMARK. The restriction, in the above statement, to infinite number of particles excludes the trivial example of a *quasi-stationary* state which is not the projection of any *steady* state, which is provided by a single point moving on the line by independent increments. In this case the state looks stationary from the perspective of the "leader": there is always just one point, at the origin. There is, however, no steady velocity $V$ such that (2.8) holds.

REMARK. Linear superpositions (of measures on the suitable $\sigma$-algebras) preserve the property of *quasi-stationarity* though not that of *steady state*—due, in the latter case, to the possible variation in the front velocities.

**3. The REM states.** We recall that for a probability measure $\rho(dx)$ on $\mathbb{R}$, a Poisson process with the density $\rho$ is a probability measure on $(\Omega, \mathcal{B})$ for which the occupation numbers for disjoint sets $A \subset \mathbb{R}$ form independent random variables, $\xi(A; \omega) \equiv \xi(A)$, with the Poisson distributions

$$\text{Prob}(\xi(A) = k) = \frac{\rho(A)^k}{k!} e^{-\rho(A)} \quad \text{and} \quad \text{mean } E(\xi(A)) = \rho(A). \tag{3.1}$$

We denote by $\mu_{s,z}(d\omega)$ the Poisson process with density $\rho_{s,z}(dx) = se^{-s(x-z)}\,dx$ on $\mathbb{R}$.

The special role of the above states in the present context is already seen in the following statement, which is based on known results.

PROPOSITION 3.1 (Based on [[3, 9]–[11]]). *For any nonlattice single-step probability distribution $g(dx)$, the collection of the steady states corresponding to the evolution by i.i.d. increments $\{h_n\}$ with the distribution $g(dh)$, as described by (2.2), consists exactly of the probability measures $\mu_{s,z}(d\omega)$ [on $(\Omega, \mathcal{B})$], with $s > 0$, $z \in \mathbb{R}$. For each of these states, the corresponding front velocity $V$ is the solution of*

$$e^{sV} = \int e^{sx} g(dx). \tag{3.2}$$

*Furthermore, with respect to $\mu_{s,z}(d\omega)$, the past increments also form an i.i.d. sequence, however with a modified distribution: conditioned on $\{\tilde{x}_n\}$, the variables $\{h_{\Pi_n}\}$ form a sequence of i.i.d. variables with the probability distribution*

$$\tilde{g}(dh) = \frac{e^{sh} g(dh)}{\int_{\mathbb{R}} e^{sy} g(dy)}. \tag{3.3}$$

Thus for these steady states the distribution of the increments changes depending on whether one looks forward or backward in time (!). In other words, the permutation $\Pi_n(\omega)$ transforms the sequence of i.i.d. variables $\{h_n\}$ into an i.i.d. sequence ($\{h_{\Pi_n}\}$) with a different distribution. (Of course this is possible only in infinite systems.)

PROOF OF PROPOSITION 3.1. The evolution by independent increments is well known, and easily seen, to take a Poisson point process into another such process with the density modified through convolution ($\rho \mapsto \rho * g$). Therefore, just the steady-state property of the states $\mu_{s,z}$ is an elementary consequence of the behavior of the exponential density under convolutions. However, for the more complete statement made above it is useful to appreciate the following observation, concerning two possible ways of viewing the collection of variables $\widetilde{\omega} = \{(x_n, h_n)\}$. The following are equivalent constructions of a point process in $\mathbb{R} \times \mathbb{R}$:

(i) A collection of points $\{x_n\}$ is generated via a Poisson process on $\mathbb{R}$, with the density $\rho(dx)$, and then to each point is attached, by its order on $\mathbb{R}$, a random variable $\{h_n\}$, taken from an i.i.d. sequence with the distribution $g(dh)$.

(ii) The configuration is generated directly as a Poisson point process in $\mathbb{R} \times \mathbb{R}$, with the two-dimensional density $\rho(dx)g(dh)$.

The transition of the perspective from (ii) to (i) requires only the second factor in the product measure on $\mathbb{R} \times \mathbb{R}$ to be normalized $\int_{\mathbb{R}} g(dh) = 1$.

Now, the map $(x, h) \mapsto (x + h, h) \equiv (\tilde{x}, h)$ takes the Poisson process describing $\widetilde{\omega}$ into another Poisson process on $\mathbb{R} \times \mathbb{R}$, which yields the joint distribution of the "new" positions paired with the steps "just taken." In case of $\mu_{s,z}(dx) \times g(dh)$, the density of the new process is: $se^{-sx}\,dx g(dh) = se^{-s(\tilde{x}-h)}\,d\tilde{x} g(dh)$. This can also be written as a product $[\int e^{sy} g(dy)] se^{-s\tilde{x}}\,d\tilde{x} \times \frac{e^{sh} g(dh)}{\int e^{sy} g(dy)}$, where now the second factor is properly normalized. By the previous observation it immediately follows that:

(i) The positions after the jump $\{\tilde{x}_n\}$ are distributed as a Poisson process on $\mathbb{R}$ with the modified density $\tilde{\xi}(dx) = [\int e^{sy} g(dy)] se^{-s\tilde{x}}\,d\tilde{x} = se^{-s(\tilde{x}-V)}\,d\tilde{x}$, that is, $\{\tilde{x}_n\}$ have the same distribution as $\{x_n\} + V$ with $V$ satisfying (3.2).

(ii) When conditioned on the configuration $\{\tilde{x}_n\}$, the jumps just taken are generated by an independent process on $\mathbb{R}$ with the probability density given by (3.3), as claimed.

For the converse statement, that is, to prove that all steady states are of the REM type, one may first note that if $\mu(d\omega)$ is a *steady state* for the dynamics corresponding to $g(dx)$, with the front velocity $V$, then $\mu$ is *stationary* under the dynamics corresponding to a shifted single-step distribution: $g_V(dx) = g(d(x - V))$. The classification of stationary states, and hence also steady states, is found in [9], where it is implied that any stationary measure is a superposition of Poisson processes whose random density solves the equation $\rho = \rho * g$. As established in [3], for nonlattice $g(dh)$ the only solutions in the space of positive measures are measures of the form $\xi(dh) = [Ae^{-sh} + B]\,dh$. The condition that the typical configurations be bounded on the positive side imply that $s > 0$ and $B = 0$. □

Having introduced the REM states, we are ready to formulate the main result.

## 4. Classification of quasi-stationary states.

DEFINITION 4.1. A probability measure $\mu$ on $\Omega$ is $g$-regular if for almost every $T\omega = \{\omega, \{h_n\}_n\}$, with respect to $\mu(d\omega)\prod_{n\in\mathbb{Z}} g(dh_n)$, the point configuration $\{x_n + h_n\}_n$ is locally finite, with a finite maximum.

The $g$-regularity of $\mu$ means that with probability 1 the configuration obtained through the independent increments has a maximum and can be ordered. This is a preliminary requirement for the possible quasi-stationarity of $\mu$. It is easy to see that a sufficient condition for $g$-regularity is met in the situation discussed next. The general sufficient condition is the finiteness, for all $x \in \mathbb{R}$, of

$$(4.1) \qquad \mathbb{E}(\text{card}\{n : x_n + h_n \geq x\}) = \mathbb{E}_\mu\left(\sum_n \text{Prob}_g(h_n \geq x - x_n)\right).$$

In the following, to simplify the exposition and to avoid confusion we will always assume that at $t=0$ the rightmost particle in the configuration is at $x_1=0$ (we can do this without loss of generality); we will denote by $x_n$ the positions of the particles at $t=0$ and by $y_n$ the positions of the particles at $t=\tau$.

Following is our main result:

THEOREM 4.2. *Let $g$ be a probability measure with a density on $\mathbb{R}$ and let $\mu$ be a probability measure on $(\Omega, \widetilde{\mathcal{B}})$, satisfying*

$$\int e^{sx} g(x)\, dx < \infty \qquad \forall\, s \in \mathbb{R} \tag{4.2}$$

*and*

$$\begin{aligned} &\mathbb{E}_\mu(\{\sharp \text{ of particles within} \\ &\qquad \text{distance } y \text{ of the leading particle}\}) \le Ae^{\lambda y} \qquad \forall\, y \ge 0 \end{aligned} \tag{4.3}$$

*for some $\lambda > 0$ and $A < \infty$. If $\mu$ is quasi-stationary with respect to the dynamics corresponding to independent increments with the distribution $g$, then it is supported on Poisson processes with densities $se^{-sx}\,dx$, $s > 0$.*

The meaning of the theorem is that the probability space $\Omega$ can be split into pieces and the process on each piece of $\Omega$ is a Poisson process with a density $se^{-sx}\,dx$ for a particular $s$.

In the proof we shall use the fact that point processes are uniquely determined by their probability generating functionals (as discussed in [4]). Our derivation of Theorem 4.2 proceeds along the following steps.

1. First we note that any quasi-stationary state can be presented as the result of evolution of arbitrary duration ($\tau$) which starts from a random initial configuration, distributed by the given quasi-stationary state, and evolves through independent increments.
2. Analyzing the above dynamics, we show that for large $\tau$ the resulting distribution is asymptotic to Poisson processes with the corresponding (evolving) densities. Thus, it is shown that the quasi-stationary measure $\mu$ can be presented as the limit of a superposition of *random* Poisson processes, where the randomness is in both the Poisson measure and the resulting particle configuration (Theorem 5.1).
3. Applying a result from the theory of large deviations (Theorem A.1), and some compactness bounds which are derived from quasi-stationarity, we show that the quasi-stationary measure admits a representation as a random Poisson process, whose Poisson densities ($F$) are the Laplace transforms of (random) positive measures (Theorem 6.1). Furthermore, in this integral representation of $\mu$, $F$ may be replaced by its convolution with $g$, followed by a normalizing shift.

4. For the last step we develop some monotonicity tools (Section 7), for which the underlying fact is that under the convolution dynamics the Laplace measures increase their relative concentration on the higher values of the Laplace parameter (Theorem 7.3). This corresponds to the statement that unless the function $F$ is a strict exponential, under the convolution dynamics the function $F$ becomes steeper, and the distribution of the gaps is shifted down. Using a strict monotonicity argument, we show that quasi-stationarity requires the measures in the above superposition to be supported on pure exponential functions (or, alternatively stated, functions whose Laplace measure is concentrated at a point).

The final implication is that the quasi-stationary measure is a superposition of REM measures, as asserted in Theorem 4.2.

Let us remark that Section 7 may be of independent interest. It is noted there that within the space of decreasing functions which are the Laplace transforms of positive measures on $[0,\infty)$, convolution with a probability measure makes any function steeper, in the sense presented below, except for the pure exponentials on which the effect of such a convolution is only a shift.

## 5. Representation of $\mu$ as a random Poisson process.

5.1. *"Poissonization"—the statement.* Let $\mathcal{F}$ be the space of monotone decreasing, continuous functions $F:\mathbb{R}\to[0,\infty]$, with $F(x)\to 0$ for $x\to\infty$ and $F(x)\to\infty$ for $x\to-\infty$. We regard a function $F\in\mathcal{F}$ as normalized if

$$F(0)=1, \tag{5.1}$$

and denote by $\mathcal{N}$ the normalizing shift: $\mathcal{N}:F(\cdot)\mapsto F(\cdot+z_F)$, with $z_F=\sup\{z\in\mathbb{R}:F(z)\geq 1\}$.

For each $F\in\mathcal{F}$, the Poisson process on $\mathbb{R}$ which corresponds to the measure $(-)dF$ will almost surely exhibit a configuration which can be ranked in the decreasing order of $\mathbb{R}$. The probability that there is no particle above $x\in\mathbb{R}$ is $\exp(-F(x))$. Conditioned on the location of the leading particle $(x)$, the rest are distributed by a Poisson process on $(-\infty,x]$ with the density $d(-F)$. Thus, the MPGFL [defined in (2.6)] of the Poisson process with density $F$, which we shall denote by $\widehat{G}_F(f)$, is given by

$$\widehat{G}_F(f)=\int_{-\infty}^{\infty} d[e^{-F(x)}]\exp\left\{-\int_{-\infty}^{x}(1-e^{-f(x-y)})\,d(-F(y))\right\}. \tag{5.2}$$

Let us note that

$$\widehat{G}_F(f)=\widehat{G}_{\mathcal{N}F}(f), \tag{5.3}$$

since the probability distribution of the gaps is not affected by uniform shifts.

For the purpose of the following theorem let $S_\tau$ be a random variable with the probability distribution $P(S_\tau \geq y) = \int I[\sum y_j \geq y] g(y_1) \cdots g(y_\tau)\, dy_1 \cdots dy_\tau$. We associate with each configuration $\omega$, and $\tau \in \mathbb{N}$, the function

$$F_{\omega;\tau}(x) = \sum_m P(S_\tau \geq x - x_m), \tag{5.4}$$

and denote by $z_{\omega,\tau}$ the position at which

$$F_{\omega;\tau}(z_{\omega,\tau}) = 1. \tag{5.5}$$

One may note that $F_{\omega;\tau}(x)$ is the expected number of particles on $[x, \infty)$ for the configuration which will be obtained from $\omega$ after $\tau$ steps of evolution with independent increments. If the support of $g(y)$ is not bounded, one may easily find configurations for which $F_{\omega;\tau}(\cdot)$ diverges. However, if the measure $\mu$ is $g$-regular, then a.s. $F_{\omega;\tau}(\cdot) < \infty$. Furthermore, we shall see that if $\mu$ is quasi-stationary, then the position of the front after $\tau$ steps can be predicted to be in the vicinity of $z_{\omega;\tau}$—up to a fluctuation whose distribution remains stochastically bounded (i.e., forms a "tight" sequence) as $\tau \to \infty$.

The main result of this section is:

Theorem 5.1. *Let $\mu$ be a $g$-regular quasi-stationary measure, for the independent evolution by steps with some common probability distribution which has a density $g(u)$. Then for every positive function $f$ of compact support in $\mathbb{R}$,*

$$\begin{aligned} \widetilde{G}_\mu(f) &= \lim_{\tau \to \infty} \int_\Omega \mu(d\omega) \widehat{G}_{\mathcal{N}F_{\omega;\tau}}(f) \\ &= \lim_{\tau \to \infty} \int_\Omega \mu(d\omega) \widehat{G}_{g * \mathcal{N}F_{\omega;\tau}}(f) \end{aligned} \tag{5.6}$$

*where $\widetilde{G}_\mu(f)$ is the modified probability generating functional defined in* (2.6).

This statement implies that the measure $\mu$ is, in the "weak sense," a limit of random Poisson processes, of measures corresponding to the random functions $\mathcal{N}F_{\omega;\tau}(\cdot)$ whose probability distribution is induced from $\mu$ through their dependence on $\omega$.

Let us note that this result is related to—but not covered by—the known statement that any limit of a sequence of point processes which is derived through successive random independent increments is a mixed Poisson process (e.g., [4], Theorem 9.4.2). Unlike in that case, the time evolution considered here incorporates shifts according to the position of the leading point (and the limiting process is not stationary under translations).

The rest of this section is devoted to the proof of this assertion, for which we need some preparatory estimates.

First let us make the following observation:

Lemma 5.2. *Any quasi-stationary measure is supported on configurations with either exactly one particle, or infinitely many.*

Proof. The statement is a simple consequence of the spreading of the probability distribution of the sum of independent increments, that is, of the variable $S_\tau$. For example, one may consider the function

$$Y_\mu^{(2)}(y) = \mu(\{y_1 - y_2 \geq y\}). \tag{5.7}$$

By the dominated-convergence theorem, $Y_\mu^{(2)}(y) \underset{y\to\infty}{\longrightarrow} 0$. However, for any finite number of particles, the probability that after $\tau$ steps the smallest gap will exceed $y$ tends to 1 as $\tau \to \infty$. Thus finite configurations of more than one particle can carry only zero probability in any quasi-stationary measure. Of course, a measure with exactly one particle is quasi-stationary. $\square$

5.2. *Some auxiliary estimates.* Given an initial configuration $\omega = \{x_n\}$, the probability distribution of the position of the leading particle after $\tau > 0$ steps is $dP_\omega^{(\tau)}(x)$, with

$$\begin{aligned} P_\omega^{(\tau)}(x) &= \text{Prob}(\{\text{at time } \tau \text{ all particles are on } (-\infty, x]\}) \\ &= \prod_n [1 - P(S_\tau \geq x - x_n)]. \end{aligned} \tag{5.8}$$

We shall need to compare $dP_\omega^{(\tau)}(x)$ with the probability distribution associated with the function

$$\widetilde{P}_\omega^{(\tau)}(x) = \exp\left\{-\sum_n P(S_\tau \geq x - x_n)\right\} = e^{-F_{\omega;\tau}(x)}. \tag{5.9}$$

Remark. It is instructive to note that $d\widetilde{P}_\omega^{(\tau)}(x)$ is the probability distribution of the maximum of a modified process, in which at first each particle is replaced by a random number of descendents, with the Poisson distribution $p_n = e^{-1}/n!$, and then each particle evolves by $\tau$ independent increments, as in the Indy-500 model. Conditioned on the starting configuration, the modified process is (instantaneously) a Poisson process. The probability that its maximum is in $(-\infty, x]$ is given by

$$\begin{aligned} &\prod_n \left[\sum_n \frac{e^{-1}}{n!}(1 - P(S_\tau \geq x - x_n))^n\right] \\ &\quad = \exp\left\{-\sum_n P(S_\tau \geq x - x_n)\right\} = \widetilde{P}_\omega^{(\tau)}(x). \end{aligned} \tag{5.10}$$

Our first goal is to show that the probability measures $dP_\omega^{(\tau)}(x)$ and $d\widetilde{P}_\omega^{(\tau)}(x)$ are "typically"—in a suitable stochastic sense—asymptotic to each other as $\tau \to \infty$. This statement is not true for some $\omega$, and it is not difficult to construct examples of configurations for which it does not hold. We note that it is easy to show that the step described by the graph of $P_\omega^{(\tau)}(\cdot)$ remains tight, in the sense that the width of the intervals $\{x : \delta \le P_\omega^{(\tau)}(x) \le 1-\delta\}$ does not spread indefinitely, as $\tau \to \infty$.

LEMMA 5.3. *For any quasi-stationary measure $\mu$:*

$$\mathbb{E}_\mu\left(\int_{-\infty}^{\infty} \sup_n P(S_\tau \ge x - x_n)\, dP_\omega^{(\tau)}(x)\right) \underset{\tau\to\infty}{\longrightarrow} 0. \tag{5.11}$$

*Furthermore,*

$$\mathbb{E}_\mu\left(\sup_x |\widetilde{P}_\omega^{(\tau)}(x) - P_\omega^{(\tau)}(x)|\right) \underset{\tau\to\infty}{\longrightarrow} 0. \tag{5.12}$$

REMARK. The supremum in (5.11) is clearly attained at $n=1$ (by monotonicity). Since $dP_\omega^{(\tau)}(x)$ is a probability measure, and the c.d.f. of $S_\tau$ is a bounded function, the statement means that the maximum typically occurs in a region whose a priori probability of being reached by any specific point is asymptotically zero.

PROOF OF LEMMA 5.3. Due to the spreading property of convolutions of probability measures (see [4], Lemma 9.4.1), for any $D < \infty$

$$b(\tau, D) = \sup_x P(x \le S_\tau < x + D) \underset{\tau\to\infty}{\longrightarrow} 0.$$

Observe that $P_\omega^{(\tau)}(x) \le \widetilde{P}_\omega^{(\tau)}(x) \le 1$ for all $x$. Let us pick $\lambda > 0$ such that

$$e^{-x(1+\lambda x)} \le 1 - x \qquad \forall\, x \in [0, \tfrac{1}{2}].$$

Thus if $P(S_\tau \ge x) \le \frac{1}{2}$, we have

$$\widetilde{P}_\omega^{(\tau)}(x) \le P_\omega^{(\tau)}(x)^{1/[1+\lambda P(S_\tau \ge x)]}. \tag{5.13}$$

Suppose that $x$ is such that $P(S_\tau \ge x) \le \varepsilon$. Then

$$\widetilde{P}_\omega^{(\tau)}(x) - P_\omega^{(\tau)}(x) \le \sup_{u\in[0,1]} |u^{1/(1+\lambda\varepsilon)} - u|.$$

Suppose that $x$ is such that $P(S_\tau \ge x - x_1) \ge \varepsilon$.

Let $n_0 = \frac{2}{\varepsilon} \ln \frac{1}{\varepsilon}$. Then for all $\tau$ large enough and for all $n \le n_0$,

$$b(\tau, -x_n) \le \frac{\varepsilon}{2}.$$

Consequently,

$$P(S_\tau \geq x - x_n) \geq P(S_\tau \geq x) - b(\tau, x_n) \geq \frac{\varepsilon}{2}.$$

Then

$$-\sum_n P(S_\tau \geq x - x_n) \leq -\sum_{n=0}^{(2/\varepsilon)\ln(1/\varepsilon)} \frac{\varepsilon}{2} \leq -\ln\frac{1}{\varepsilon},$$

and therefore

$$\widetilde{P}_\omega^{(\tau)}(x) \leq e^{-\ln(1/\varepsilon)} \leq \varepsilon. \tag{5.14}$$

So in this case we obtain

$$P_\omega^{(\tau)}(x) \leq \varepsilon.$$

□

Putting the above together, we have:

Lemma 5.4. *If $\mu$ be a quasi-stationary measure, then for each $\varepsilon > 0$,*

$$\mu(\{\omega : \operatorname{dist}(dP_\omega^{(\tau)}, d\widetilde{P}_\omega^{(\tau)}) \geq \varepsilon\}) \underset{\tau\to\infty}{\longrightarrow} 0, \tag{5.15}$$

*where* dist *is the distance between the two measures, defined as*

$$\operatorname{dist}(dP, d\widetilde{P}) = \sup_h \left\{ \left| \int h(x)\, dP(x) - \int h(x)\, d\widetilde{P}(x) \right| \Big/ \|h\|_\infty \right\}. \tag{5.16}$$

Proof. The distributions $dP_\omega^{(\tau)}(x)$ and $d\widetilde{P}_\omega^{(\tau)}(x)$ can be written as

$$\begin{aligned} dP_\omega^{(\tau)}(x) &= \sum_k \frac{dP(S_\tau \geq x - x_k)}{1 - P(S_\tau \geq x - x_k)} \prod_n [1 - P(S_\tau \geq x - x_n)], \\ d\widetilde{P}_\omega^{(\tau)}(x) &= \sum_k dP(S_\tau \geq x - x_k) \times \exp\left\{ -\sum_n P(S_\tau \geq x - x_n) \right\}. \end{aligned} \tag{5.17}$$

By Lemma 5.3 we obtain that

$$|P_\omega^{(\tau)}(x) - \widetilde{P}_\omega^{(\tau)}(x)| \leq \varepsilon \qquad \forall x.$$

If $P(S_\tau \geq x) \leq \varepsilon$, then we obtain by the same arguments as in the previous lemma that

$$|d\widetilde{P}_\omega^{(\tau)}(x) - dP_\omega^{(\tau)}(x)| \leq \varepsilon\, d\widetilde{P}_\omega^{(\tau)}(x).$$

Integrating with respect to $\frac{h}{\|h\|_\infty}$ over the $x$ such that $P(S_\tau \geq x) \leq \varepsilon$, we obtain that the result is small.

If $P(S_\tau \geq x) > \varepsilon$, then

$$\prod_{n\neq k}[1 - P(S_\tau \geq x - x_n)] \leq \varepsilon \quad \text{and} \quad \exp\left\{-\sum_n P(S_\tau \geq x - x_n)\right\} \leq \varepsilon.$$

Consequently for such $x$

$$\begin{aligned}\int \frac{h(x)}{\|h\|_\infty}\, d\widetilde{P}_\omega^{(\tau)}(x) &\\ &\leq \sqrt{\varepsilon}\int \sum_n dP(S_\tau \geq x - x_n)\exp\left\{-\frac{1}{2}\sum_n P(S_\tau \geq x - x_n)\right\}\\ &\leq \text{const}\,\sqrt{\varepsilon},\end{aligned}$$

and also using $1 - x \leq e^{-x}$ for $x > 0$

$$\begin{aligned}\int \frac{h(x)}{\|h\|_\infty}\, dP_\omega^{(\tau)}(x) &\\ &\leq \sqrt{\varepsilon}\int d\sum_n P(S_\tau \geq x - x_n)\exp\left\{-\frac{1}{2}\sum_n P(S_\tau \geq x - x_n)\right\}\\ &\leq \text{const}\,\sqrt{\varepsilon}.\end{aligned}$$

□

5.3. "*Poissonization*"*—the proof.* We are now ready to prove the main result of this section.

Proof of Theorem 5.1. Due to the quasi-stationarity of the measure $\mu$, one may evaluate $\widetilde{G}_\mu(f)$ by taking the average of the future expectation value of $\exp\{-\sum_n f(y_1 - y_n)\}$, corresponding to the configuration $\omega$ as it appears at time $t = 0$.

In the following argument we fix the (nonnegative) "test function" $f$, and take $D < \infty$ such that $\operatorname{supp} f \subset [-D, 0]$. In the approximations which follow we use the fact that $\exp\{-\sum f(y_1 - y_n)\}$ is a bounded function ($\leq 1$), which is integrated against a probability measure. As before, $\omega$-dependent quantities are denoted $o(1)$ if in the limit $\tau \to 0$ they tend to 0 "in law," that is, the probability distribution which they inherited from $\omega$ is nonzero only for $[0, \varepsilon]$ for any $\varepsilon > 0$.

The conditional expectation of the future value of $\exp\{-\sum f(y_1 - y_n)\}$, conditioned on the initial configuration $\omega$, is

$$\begin{aligned}&\mathbb{E}_\omega\Big(\exp\Big\{-\sum f(y_1 - y_n)\Big\}\Big)\\ &\qquad = \int_{-\infty}^{\infty} e^{-f(0)}\, dP(S_\tau \geq x - x_k)\end{aligned} \tag{5.18}$$

$$\times \prod_{n\neq k}[1-P(S_\tau \geq x-x_n)]$$

$$\times \prod_{n\neq k} \frac{\int_{-\infty}^{x} e^{-f(x-y)}\,dP(S_\tau \geq y-x_n)}{(1-P(S_\tau \geq x-x_n))},$$

where $dP(S_\tau \geq x-x_k)$ is the probability that the $k$th particle is at $x$ at time $\tau$, $\prod_{n\neq k}[1-P(S_\tau \geq x-x_n)]$ is the probability that other particles are at $(-\infty, x]$ at time $\tau$, and $\frac{\int_{-\infty}^{x} e^{-f(x-y)}\,dP(S_\tau \geq y-x_n)}{(1-P(S_\tau \geq x-x_n))}$ is the expectation of $e^{-f(x-y_n)}$ given that the particle which is at $x_n$ at $t=0$ is at $(-\infty, x]$ at time $\tau$.

As in the previous discussion, the contribution of $x$ such that $P(S_\tau \geq x) \geq \varepsilon$ to the integral in (5.18) is negligible.

Consider $x$ such that $P(S_\tau \geq x) \leq \varepsilon$. We can write

$$\prod_{n\neq k} \frac{\int_{-\infty}^{x} e^{-f(x-y)}\,dP(S_\tau \geq y-x_n)}{1-P(S_\tau \geq x-x_n)}$$

$$= \prod_{n\neq k}\left[1-\frac{\int_{-\infty}^{x}(1-e^{-f(x-y)})\,dP(S_\tau \geq y-x_n)}{1-P(S_\tau \geq x-x_n)}\right] \tag{5.19}$$

$$= (1+o(1))\exp\left\{-\int_{-\infty}^{x}(1-e^{-f(x-y)})\,d\left(\sum_n P(S_\tau \geq x-x_n)\right)\right\}.$$

As noted in (5.3), the normalizing shift has no effect on $\widehat{G}_F(f)$. The result is the first of the two equations in (5.6). The second equation is an immediate corollary of the first one, since

$$g * F_{\omega;\tau} = F_{\omega;\tau+1}. \tag{5.20}$$ □

For a later use, let us note that the arguments used in the above discussion readily imply the following two bounds.

COROLLARY 5.5. *For any $\varepsilon > 0$, there is $W(\varepsilon) < \infty$ such that*

$$\mathbb{E}_\mu\left(\int_{|x|>W(\varepsilon)} d[e^{-\mathcal{N}F_{\omega;\tau}(x)}]\right) \leq \varepsilon \tag{5.21}$$

*and*

$$\mathbb{E}_\mu\left(\int_{|x|>W(\varepsilon)} d[e^{-g*\mathcal{N}F_{\omega;\tau}(x)}]\right) \leq \varepsilon. \tag{5.22}$$

PROOF. Let $f = \mathrm{I}_{[0,W(\varepsilon)]}$. Denote

$$\phi(W(\varepsilon)) = \mathbb{E}_\mu[e^{-\mathrm{I}_{[0,W(\varepsilon)]}(y_1-y_n)}]. \tag{5.23}$$

Since $\mathrm{I}_{[0,W(\varepsilon)]}(x) \underset{W(\varepsilon)\to\infty}{\longrightarrow} 1$ for $x \in \mathbb{R}$ and since, in a typical configuration, the number of particles within distance $W(\varepsilon)$ behind the leader increases to $\infty$ as $W(\varepsilon)$ increases, $\phi(W(\varepsilon))$ must decay monotonically to 0 as $W(\varepsilon)$ increases. By taking $f = \mathrm{I}_{[0,W(\varepsilon)]}$, we see that

$$\begin{aligned}\phi(W(\varepsilon)) = \int \mu(d\omega) \int_{-\infty}^{\infty} de^{-\mathcal{N}F_{\omega;\tau}(x)} \\ \times e^{-(1-e^{-1})(\mathcal{N}F_{\omega;\tau}(x-W(\varepsilon))-\mathcal{N}F_{\omega;\tau}(x))} + O(\varepsilon_\tau).\end{aligned} \tag{5.24}$$

We can get an estimate on $\mathcal{N}F_{\omega;\tau}(W(\varepsilon))$ from (5.24) by restricting the range of integration from $W(\varepsilon)$ to $\infty$ and using that $\mathcal{N}F_{\omega;\tau}(x - W(\varepsilon)) - \mathcal{N}F_{\omega;\tau}(x) \leq 1$. Then, for $x \geq W(\varepsilon)$ we obtain

$$\begin{aligned}\phi(W(\varepsilon)) &\geq \int \mu(d\omega) \int_{W(\varepsilon)}^{\infty} de^{-\mathcal{N}F_{\omega;\tau}(x)} \\ &\quad \times e^{-(1-e^{-1})(\mathcal{N}F_{\omega;\tau}(x-W(\varepsilon))-\mathcal{N}F_{\omega;\tau}(x))} + O(\varepsilon_\tau) \\ &\geq e^{-(1-e^{-1})} \int \mu(d\omega)(1 - e^{-\mathcal{N}F_{\omega;\tau}(W(\varepsilon))}) + O(\varepsilon_\tau).\end{aligned} \tag{5.25}$$

Similarly, by restricting the range of integration from 0 to $\infty$ and using that $\mathcal{N}F_{\omega;\tau}(x - W(\varepsilon)) - \mathcal{N}F_{\omega;\tau}(x) \leq \mathcal{N}F_{\omega;\tau}(-W(\varepsilon))$ for $x \geq 0$, we obtain

$$\begin{aligned}\phi(W(\varepsilon)) &\geq \int \mu(d\omega) \int_{0}^{\infty} de^{-\mathcal{N}F_{\omega;\tau}(x)} \\ &\quad \times e^{-(1-e^{-1})(\mathcal{N}F_{\omega;\tau}(x-W(\varepsilon))-\mathcal{N}F_{\omega;\tau}(x))} + O(\varepsilon_\tau) \\ &\geq (1 - e^{-1}) \int \mu(d\omega) e^{-(1-e^{-1})\mathcal{N}F_{\omega;\tau}(-W(\varepsilon))} + O(\varepsilon_\tau).\end{aligned} \tag{5.26}$$

Equations (5.25) and (5.26) prove the first part of the corollary.

To prove (5.22) we observe that from the previous part it follows that for all $\tau$ large enough, and sufficiently large $W(\varepsilon)$,

$$\mathbb{E}_\mu \int_{|x| \geq W(\varepsilon)/2} de^{-\mathcal{N}F_{\omega;\tau+1}}(x) \leq \frac{\varepsilon}{2}.$$

Since for sufficiently large $W(\varepsilon)$ and $\omega$ in a set of measure $1 - \frac{\varepsilon}{2}$,

$$z_{\omega,\tau+1} - z_{\omega,\tau} \leq \frac{W(\varepsilon)}{2},$$

we obtain that

$$\mathbb{E}_\mu \int_{|x| \geq W(\varepsilon)} de^{-g * \mathcal{N}F_{\omega;\tau}(x)} \leq \varepsilon.$$

□

Corollary 5.5 will be used for an approximation of $\widehat{G}_{\mathcal{N}F_{\omega;\tau}}(f)$ by a quantity which has better continuity properties as a functional of $F$.

**6. The Poisson density as a Laplace transform of a random positive measure.** We shall next show that the quasi-stationary measure $\mu$ can be presented as equivalent to a random Poisson process whose density is the Laplace transform of a random positive measure on $\mathbb{R}$. [Due to the invariance of $\widetilde{\mathcal{B}}$ under uniform shifts, with no additional restriction the measures may be adjusted so that $\rho(\mathbb{R}) = 1$.]

Let $\mathcal{M}$ be the space of finite measures on $[0, \infty)$. To each $\rho \in \mathcal{M}$ we associate the Laplace transform function

$$R_\rho(x) = \int_0^\infty e^{-xu} \rho(du). \tag{6.1}$$

We denote by $\mathcal{F}_L$ the space of such functions, that is, $\mathcal{F}_L = \{R_\rho(\cdot) | \rho \in \mathcal{M}\}$.

We shall need to consider "ensemble averages" over randomly chosen elements of $\mathcal{M}$. These are described by probability measures on $\mathcal{M}$, which would always be understood to be defined on the natural $\sigma$-algebra on $\mathcal{M}$, for which the measures of intervals, $\rho([a, b])$, are measurable functions of $\rho$. Our goal in this section is to prove the following statement.

Theorem 6.1. *Under the assumptions of Theorem 4.2, there exists a probability measure, $\nu(d\rho)$, on $\mathcal{M}$ such that for any compactly supported positive function $f$ on $\mathbb{R}$,*

$$\widetilde{G}_\mu(f) = \int_{\mathcal{M}} \nu(d\rho) \widehat{G}_{R_\rho}(f), \tag{6.2}$$

*and furthermore,*

$$\widetilde{G}_\mu(f) = \int_{\mathcal{M}} \nu(d\rho) \widehat{G}_{R_\rho * g}(f). \tag{6.3}$$

For Laplace transform functions $F = R_\rho$, shifts correspond to transformations of the form

$$\rho(du) \Longrightarrow e^{-\alpha u} \rho(du), \tag{6.4}$$

and the normalization condition (5.1) corresponds to $\rho(\mathbb{R}) = 1$, that is, $\rho \in \mathcal{M}$ being a *probability measure*. In view of the invariance (5.3), this normalization condition may be freely added as a restriction of the support of $\nu(d\rho)$ in the statement of Theorem 6.1.

While the result presented in the previous section required only quasi-stationarity, we shall now make use of the additional assumptions listed in the main theorem (Theorem 4.2).

In the derivation of Theorem 6.1 we shall apply what may be regarded as the principle of the equivalence of ensembles, in the language of statistical mechanics. Specifically, we need the following result, which, as is explained in the Appendix, is a refinement of the "Bahadur–Rao theorem" of large deviation theory.

THEOREM A.1. *Let $u_1, u_2, \ldots$ be i.i.d. random variables with expectation $E_g u$ and a common probability distribution $g(u)$, which has a density and a finite moment generating function, $\int e^{\eta u} g(u)\, du \equiv e^{\Lambda(\eta)} < \infty$ for all $\eta$. Then, for any $0 < K < \Lambda'(\infty)$ and $0 < \beta < \frac{1}{2}$ there is $\varepsilon_{\tau;K,\beta} \underset{\tau\to\infty}{\longrightarrow} 0$ such that for all $q \in [E_g u, K]$ and $|x| \le \tau^\beta$,*

$$\frac{\mathrm{Prob}(\{u_1 + u_2 + \cdots + u_\tau \ge x + q\tau\})}{\mathrm{Prob}(\{u_1 + u_2 + \cdots + u_\tau \ge q\tau\})} = e^{-\eta x}[1 + O(\varepsilon_{\tau;K,\beta})], \tag{6.5}$$

*with $\eta = \eta(q)$ determined by the condition*

$$q = \frac{\int u e^{\eta u} g(u)\, du}{\int e^{\eta y} g(y)\, dy}. \tag{6.6}$$

*In our analysis we shall need a bound on the front velocity, and on the possible propagation of particles from the far tail.*

LEMMA 6.2. *Let $\mu$ be a quasi-stationary $g$-regular measure with a density satisfying the assumptions (4.2) and (4.3) of Theorem 4.2. Then:*

(i) *For any $\tau$ large enough, for $\omega$ in a set of measure $1 - \varepsilon$,*

$$z_{\omega;\tau} \le \frac{S}{2\lambda}\tau + \mathrm{const} \qquad \text{where } S = \ln \int e^{2\lambda x} g(x)\, dx. \tag{6.7}$$

(ii) *There exist $\alpha_\mu(M)$ and $\beta_g(\tau)$ such that the probability of the complement of the event*

$$A_{\tau;D,K,M} = \{\omega : \text{the configuration obtained after } \tau \text{ steps will have not more than } M \text{ particles with } y_n \ge y_1 - D, \text{ and all of them made a total jump less than } K\tau + z_{\omega,\tau} - x_n \text{ in time from } 0 \text{ to } \tau\}$$

*satisfies*

$$\mathrm{Prob}(A^c_{\tau;D,K,M}) \le \alpha_\mu(M, D) + \beta_g(\tau) + C_{g,\mu} e^{-\delta(K-K_0)\tau} \tag{6.8}$$

*with $\alpha_\mu(M, D) \underset{M\to\infty}{\longrightarrow} 0$ for each $D < \infty$, $\beta_g(\tau) \underset{\tau\to\infty}{\longrightarrow} 0$ and $\delta > 0$.*

REMARK. In the proof below we shall apply the last bound in the double limit: $\lim_{K\to\infty} \lim_{\tau\to\infty}$, with $M$ chosen so that $1 \ll M \ll K$.

Proof of Lemma 6.2. (i) By (4.3) and Markov inequality,

$$P_\mu(\{\sharp(-x_n)\le m\}\ge e^{2\lambda m})\le e^{-\lambda m}.$$

Therefore by the Borel–Cantelli lemma,

$$P_\mu\{\{\sharp(-x_n)\le m\}> e^{2\lambda m}\text{ i.o.}\}=0.$$

This implies that for any $\varepsilon$ there exists $m_0$ such that on a set of $\omega$ of measure $1-\varepsilon$, $\{\sharp(-x_n)\le m\}\le e^{2\lambda m}$ for all $m\ge m_0$.

Using the definition of $F_{\omega;0}$ we obtain

$$F_{\omega;0}(x)\le e^{-2\lambda\min(x,-m_0)}\qquad\forall\, x<0. \tag{6.9}$$

Therefore,

$$\begin{aligned}F_{\omega;\tau}(x)&\le F_{\omega;0}*g^{(*\tau)}(x)\\&\le \text{const}\int e^{-2\lambda(x-y)}g^{(*\tau)}(y)\,dy\le \text{const}\,e^{-2\lambda x+\tau S}.\end{aligned} \tag{6.10}$$

For $x=\frac{S\tau}{2\lambda}+\text{const}$ we thus obtain $F_{\omega;\tau}(x)\le 1$. It follows by definition that $z_{\omega,\tau}\le\frac{S\tau}{2\lambda}+\text{const}$.

(ii) The probability that the first condition does not hold in the definition of $A_{\tau;D,K,M}$ is, by the quasi-stationarity of $\mu$,

$$\begin{aligned}\alpha_\mu(M,D)=\mu(\omega:&\text{more than } M \text{ particles}\\&\text{are within distance } D \text{ of the leader at } t=0).\end{aligned} \tag{6.11}$$

This quantity vanishes for $M\to\infty$ because the number of particles in $[y_1-D,y_1]$ is almost surely finite.

To estimate the remaining probability of the complement of the event $A_{\tau;D,K,M}$ we split it into two cases, based on the distance which the front advances in time $\tau$. That distance is at least the total displacement of the particle which is initially at 0. The probability that this displacement is less than $(E_g u-1)\tau$ is dominated by the quantity

$$\text{Prob}(z_{\omega,\tau}\le(E_g u-1)\tau)\le\text{Prob}(S_\tau\le(E_g u-1)\tau)=\beta_g(\tau).$$

The choice of 1 is somewhat arbitrary, but even so, standard large deviation arguments which are applicable under the assumption (4.2) imply that $\beta_g(\tau)$ decays exponentially.

The contribution of the other case is bounded by the probability of the following event:

$$\begin{aligned}\text{Prob}(&\text{at least one of the particles of } \omega \text{ will advance in } \tau \text{ steps}\\&\text{a distance greater than } [-x_n+(E_g u-1+K)\tau])\end{aligned}$$

$$
\begin{aligned}
&\le \mathbb{E}_\mu\left(\sum_n \mathrm{Prob}(S_\tau \ge -x_n + (E_g u - 1 + K)\tau)\right) \\
&\le \mathbb{E}_\mu\left(\sum_n \mathbb{E}_g(e^{\alpha S_\tau}) e^{-\alpha[-x_n + (E_g u - 1 + K)\tau]}\right) \\
&\le \left[\int_{\mathbb{R}} e^{\alpha(u - E_g u)} g(u)\, du e^{-\alpha(K-1)}\right]^\tau \mathbb{E}_\mu\left(\sum_n e^{-\alpha[-x_n]}\right),
\end{aligned}
\tag{6.12}
$$

where $\alpha > 0$ is an adjustable constant. The last factor is finite for $0 < \alpha < \lambda$ since under the assumed exponential bound (4.3),

$$
\begin{aligned}
\mathbb{E}_\mu\left(\sum_n e^{-\alpha[-x_n]}\right) &= \mathbb{E}_\mu\left(\alpha \int dy e^{-\alpha y} \sum_n \mathrm{I}[y \ge -x_n]\right) \\
&\le \alpha \int dy e^{-\alpha y} A e^{\lambda y} = \frac{A\alpha}{\alpha - \lambda}.
\end{aligned}
\tag{6.13}
$$

The claimed estimate readily follows (choosing $\lambda < \alpha$, and defining $\delta > 0$ correspondingly). □

Proof of Theorem 6.1. Applying Theorem A.1 to the function defined by (5.4), we find that

$$
\begin{aligned}
\mathcal{N} F_{\omega;\tau}(x) &= \sum_n P(S_\tau \ge x + z_{\omega;\tau} - x_n) \\
&= \sum_{-K(\varepsilon)\tau \le x_n \le 0} P(S_\tau \ge z_{\omega;\tau} - x_n) \frac{P(S_\tau \ge x + z_{\omega;\tau} - x_n)}{P(S_\tau \ge z_{\omega;\tau} - x_n)} \\
&\quad + \sum_{x_n \le -K(\varepsilon)\tau} P(S_\tau \ge z_{\omega,\tau} - x_n + x) \\
&= \sum_{-K(\varepsilon)\tau \le x_n \le 0} P(S_\tau \ge z_{\omega;\tau} - x_n) e^{-\eta((z_{\omega;\tau} - x_n)/\tau)\cdot x}[1 + O(\varepsilon_\tau)] \\
&\quad + \sum_{x_n \le -K(\varepsilon)\tau} P(S_\tau \ge z_{\omega,\tau} - x_n + x) \\
&= \int_0^\infty \rho_{\omega;\tau}(du) e^{-ux}[1 + O(\varepsilon_\tau)] \\
&\quad + \sum_{x_n \le -K(\varepsilon)\tau} P(S_\tau \ge z_{\omega,\tau} - x_n + x),
\end{aligned}
\tag{6.14}
$$

with $\rho_{\omega;\tau}(du)$ defined as the probability measure with weights $P(S_\tau \ge z_{\omega;\tau} - x_n)$ at the points $\eta(\frac{z_{\omega;\tau} - x_n}{\tau})$.

We will now estimate the remainder term $\sum_{x_n \le -K(\varepsilon)\tau} P(S_\tau \ge z_{\omega,\tau} - x_n + x)$.

In the case when $\lim_{\eta\to\infty} \Lambda'(\eta) < \infty$ (in the case when the supremum of the support of $g(x)$ is finite), the remainder term is zero for large $K(\varepsilon)$ [e.g., if $K(\varepsilon) \ge \Lambda'(\infty)$ and $x = O(\tau^\beta)$].

In the case when $\lim_{\eta\to\infty} \Lambda'(\eta) = \infty$, the remainder term can be estimated using the large deviation arguments. By using (A.2) in Theorem A.1 and (6.9) we obtain

$$\begin{aligned}
&\sum_{x_n \le -K(\varepsilon)\tau} P(S_\tau \ge z_{\omega,\tau} - x_n + x) \\
&\qquad \le \int_{K(\varepsilon)\tau}^{\infty} P(u_1 + \cdots + u_\tau \ge y + z_{\omega,\tau} + x) e^{2\lambda y}\, dy \\
&\qquad = \int_{K(\varepsilon)\tau}^{\infty} \exp\biggl\{-\tau\Lambda^*\biggl(\frac{y + z_{\omega,\tau} + x}{\tau}\biggr)\biggr\} \\
&\qquad\qquad \times \biggl[\int_0^\infty \exp\biggl\{-\psi_\tau\biggl(\eta\biggl(\frac{y + z_{\omega,\tau} + x}{\tau}\biggr)\biggr)t\biggr\}\, dQ_\tau^{(\eta(y/\tau))}(t)\biggr] e^{2\lambda y}\, dy \qquad (6.15)\\
&\qquad = O(\varepsilon_\tau).
\end{aligned}$$

The last equality in (6.15) follows because by convexity of $\Lambda^*$

$$\Lambda^*\biggl(\frac{y + z_{\omega,\tau} + x}{\tau}\biggr) \gg \frac{2\lambda y}{\tau} \qquad \text{for all } y \ge K(\varepsilon)\tau,$$

and because the factor in the square brackets in (6.15) is small.

Therefore

$$\mathcal{N}F_{\omega;\tau}(x) = \int_0^\infty e^{-ux} \rho_{\omega;\tau}(du)(1 + O(\varepsilon_\tau)). \tag{6.16}$$

We observe that

$$\begin{aligned}
|\mathcal{N}F_{\omega;\tau}(x) - R_\rho(x)| &\le \varepsilon R_\rho(x), \\
|\mathcal{N}F'_{\omega;\tau}(x) - R'_\rho(x)| &\le \varepsilon R'_\rho(x), \\
|g * \mathcal{N}F_{\omega;\tau}(x) - g * R_\rho(x)| &\le \varepsilon g * R_\rho(x), \\
|(g * \mathcal{N}F_{\omega;\tau})'(x) - (g * R_\rho)'(x)| &\le \varepsilon (g * R_\rho)'(x).
\end{aligned} \tag{6.17}$$

Using (6.17) and Corollary 5.5 we obtain

$$\begin{aligned}
\widetilde{G}_\mu(f) &= \int_{-W(\varepsilon)}^{W(\varepsilon)} de^{-R_\rho(x)} \exp\biggl\{-\int_{-\infty}^{x} (1 - e^{-f(x-y)})(-dR_\rho(y))\biggr\} + \varepsilon \\
&= \widehat{G}_{W(\varepsilon),R_\rho}(f) + \varepsilon \qquad (6.18)
\end{aligned}$$

$$= \int_{-W(\varepsilon)}^{W(\varepsilon)} de^{-g*R_\rho(x)} \exp\left\{-\int_{-\infty}^{x} (1-e^{-f(x-y)})(-dg*R_\rho(y))\right\} + \varepsilon$$

$$= \widehat{G}_{W(\varepsilon),g*R_\rho}(f) + \varepsilon.$$

[Equation (6.18) will serve as a definition of $\widehat{G}_{W(\varepsilon),R_\rho}(f)$.]

From (6.16) we observe that for all $\omega$ in a set of measure $1-\varepsilon$ and for every $K \gg 1$ there exists a $K_1 \gg 1$ depending on $K$ such that

$$\int_{\eta(E_g u+K)}^{\infty} e^{Du} \rho_{\omega;\tau}(du) \le \sum_{x_n \le -K_1\tau} P(S_\tau \ge x + z_{\omega;\tau} - x_n). \tag{6.19}$$

We can choose $K_1$ by requiring that for all $x_n \le -K_1\tau$,

$$\frac{z_{\omega,\tau} - x_n}{\tau} \le E_g u + K,$$

for example, $K_1 = E_g u + K - \frac{S}{2\lambda}$ and we used (6.7).

From Lemma 6.2 [see (6.8)], applied with $M = \sqrt{K}$ (or any other choice with $1 \ll M \ll K$), we find that under the assumptions listed above, for any $D < \infty$ there exist $\varepsilon_D(K)$ with which

$$\begin{aligned} &\limsup_{\tau\to\infty} \mathbb{E}_\mu\left(\int_{\eta(E_g u+K)}^{\infty} e^{Du}\rho_{\omega;\tau}(du)\right) \\ &\quad \le \limsup_{\tau\to\infty} \mathbb{E}_\mu\left(\sum_{x_n \le -K_1\tau} P(S_\tau \ge z_{\omega;\tau} - x_n - D)\right) + \varepsilon_D(K) \\ &\quad \equiv \tilde{\varepsilon}_D(K) \underset{K\to\infty}{\longrightarrow} 0. \end{aligned} \tag{6.20}$$

The correspondence $\omega \mapsto \rho_{\omega;\tau}$ defines a mapping from the space of configurations $\Omega$ into the space $\mathcal{M}$, of measures on $\mathbb{R}$, with values restricted to the subset of probability measures. Corresponding to this map is one which takes the measure $\mu$ on $\Omega$ into a probability measure on $\mathcal{M}$ which we shall denote by $\nu_\tau$. By this definition, for any measurable function $\mathcal{X} : \mathcal{M} \to \mathbb{R}$,

$$\int_{\mathcal{M}} \mathcal{X}(\rho)\nu_\tau(d\rho) = \int_\Omega \mathcal{X}(\rho_{\omega;\tau})\mu(d\omega). \tag{6.21}$$

The space of probability measures on compact subsets of $\mathbb{R}$ is compact, and so is the space of probability measures on this space. While we do not have such compactness (since the measures of $\mathcal{M}$ are defined over the noncompact $\mathbb{R}$), (6.20) with any fixed $D > 0$ implies that the sequence of measures $\nu_\tau$ is *tight* and that it has a subsequence $\nu_{\tau_n}$ which converges in the corresponding "weak topology" as $\tau_n \to \infty$. Let $\nu$ be a limit of such a subsequence. [To prove the tightness of $\nu_{\tau_n}$ we observe that it is possible to show that for all $\tau$, $R_{\rho_{\omega;\tau}}(x) \le M(x)$ for some function $M(x)$ except for $\omega$ in a set of measure $\varepsilon$. The set of $\rho$ for which $R_\rho(x) \le M(x)$ is compact.] We claim that for every positive $f$ of compact support, $\operatorname{supp} f \subset [-D, 0]$,

$$\begin{aligned}
&\int_{\mathcal{M}} [\widehat{G}_{W(\varepsilon),R_\rho}(f)]\nu(d\rho) + \varepsilon \\
&\qquad = \lim_{n\to\infty} \int_{\mathcal{M}} [\widehat{G}_{W(\varepsilon),R_\rho}(f)]\nu_{\tau_n}(d\rho) + \varepsilon = \widetilde{G}_\mu(f), \\
&\int_{\mathcal{M}} [\widehat{G}_{W(\varepsilon),g*R_\rho}(f)]\nu(d\rho) + \varepsilon \\
&\qquad = \lim_{n\to\infty} \int_{\mathcal{M}} [\widehat{G}_{W(\varepsilon),g*R_\rho}(f)]\nu_{\tau_n}(d\rho) + \varepsilon = \widetilde{G}_\mu(f).
\end{aligned} \tag{6.22}$$

The weak convergence means that for any *continuous* function $\mathcal{X} : \mathcal{M} \to \mathbb{R}$,

$$\begin{aligned}
\int_{\mathcal{M}} \mathcal{X}(\rho)\nu(d\rho) &= \lim_{n\to\infty} \int_{\mathcal{M}} \mathcal{X}(\rho)\nu_{\tau_n}(d\rho) \\
&= \lim_{n\to\infty} \int_{\Omega} \mathcal{X}(\rho_{\omega;\tau_n})\mu(d\omega).
\end{aligned} \tag{6.23}$$

The continuity argument does not apply immediately to the function which we are interested in:

$$\begin{aligned}
&\widehat{G}_{W(\varepsilon),R_\rho}(f) \\
&\quad = \int_{-W(\varepsilon)}^{W(\varepsilon)} d[e^{-R_\rho(x)}] \exp\biggl\{-\int_{x-D}^{x} [1 - e^{-f(x-y)}]\, d(-R_\rho(y))\biggr\},
\end{aligned} \tag{6.24}$$

which is not continuous in $\rho$. However, $\widehat{G}_{W(\varepsilon),R_\rho}(f)$ can be approximated arbitrarily well, in the appropriate $L_1$ sense, by functionals which are continuous.

The function $\widehat{G}_{W(\varepsilon),R_\rho}(f)$ is not continuous in $\rho$. The difficulty is that $R_\rho(x)$ can be affected by small changes in the measure $\rho$ if those occur at high values of the Laplace variable $u$. However, we do obtain a continuous function by replacing $R_\rho$ in (6.24) by $R_{\mathrm{I}_K\rho}$ with

$$\mathrm{I}_K\rho(du) = \mathrm{I}_{[0,\eta(E_g u+K)]}(u) \cdot \rho(du). \tag{6.25}$$

It is easy to see that

$$\begin{aligned}
&\int_{\mathcal{M}} \biggl(\int_{\eta(E_g u+K)}^{\infty} e^{-xu}\rho(du)\biggr)\nu(d\rho) \\
&\quad \le \limsup_{\tau\to\infty} \int_{\Omega} \biggl(\int_{\eta(E_g u+K)}^{\infty} e^{-xu}\rho_{\omega;\tau}(du)\biggr)\mu(d\omega) \\
&\quad \le \tilde{\varepsilon}_x(K) \underset{K\to\infty}{\longrightarrow} 0,
\end{aligned} \tag{6.26}$$

where the first inequality is by the generalized version of Fatou's lemma, and the second is by (6.20).

Due to the fact that $f$ is compactly supported and the integration in $x$ is over $[-W(\varepsilon), W(\varepsilon)]$, the difference

$$\int_{\mathcal{M}} |\widehat{G}_{W(\varepsilon),R_\rho}(f) - \widehat{G}_{W(\varepsilon),R_{\mathrm{I}_K\rho}}(f)|\, d\nu(\rho)$$

is affected only by values of $x \in [-W(\varepsilon) - D, W(\varepsilon)]$.

Taking $x$ in this interval, we observe that (6.26) implies that

$$\int_{\eta(E_g u+K)}^{\infty} e^{-ux}\rho(du) \le \varepsilon,$$

except on the set of $\omega$ of measure $\varepsilon$.

The difference $\int_{\mathcal{M}} |\widehat{G}_{W(\varepsilon),R_\rho}(f) - \widehat{G}_{W(\varepsilon),R_{\mathrm{I}_K\rho}}(f)|\, d\nu(\rho)$ is controlled by $|R_\rho(x) - R_{\mathrm{I}_K\rho}(x)|$ and by $|R'_\rho(x) - R'_{\mathrm{I}_K\rho}(x)|$, which are small for $x \in [-W(\varepsilon) - D, W(\varepsilon)]$ except on the set of $\rho$ of measure $\varepsilon$, since $\int_{\eta(E_g u+K)}^{\infty} e^{-ux}\rho(du)$ is small.

One can verify by standard arguments that for $K$ finite $\widehat{G}_{W(\varepsilon),R_{\mathrm{I}_K\rho}}(f)$ is continuous in $\rho$, and that this continuity and the approximation bounds listed above imply (6.22), thereby proving the first part of Theorem 6.1. The second part is proved via similar arguments. □

**7. Monotonicity arguments.** In this section we develop some monotonicity tools, which will be applied to prove that if a measure $\mu$ has the properties listed in Theorem 6.1, then the corresponding measures $\rho$ are $\nu$-almost surely concentrated on points, that is, the Poisson densities $R_\rho$ are almost surely pure exponential.

7.1. *The contraction property of convolutions within* $\mathcal{F}_L$. The space $\mathcal{F}$, whose elements are positive decreasing continuous functions on $\mathbb{R}$, is partially ordered by the following relation.

DEFINITION 7.1. For $F, G \in \mathcal{F}$ we say that $G$ is *steeper* than $F$ if the level intervals of $G$ are *shorter* than those of $F$, in the sense that for any $0 \le a \le b \le \infty$,

$$(0 \le) G^{-1}(a) - G^{-1}(b) \le F^{-1}(a) - F^{-1}(b). \tag{7.1}$$

We adapt the convention that for the (monotone) functions $G \in \mathcal{F}$ the inverse is defined (for $a \ge 0$) by

$$G^{-1}(a) = \inf\{x \in \mathbb{R} : G(x) \le a\}. \tag{7.2}$$

It is easy to see that, within the class of monotone functions $\mathcal{F}$, an equivalent formulation of the relation "$G$ is steeper than $F$" is that for any $u > 0$,

$$G(x) = F(y) \quad \Longrightarrow \quad G(x+u) \le F(y+u). \tag{7.3}$$

Also equivalent is such a principle with the reversed inequality and $u < 0$.

Of particular interest for us is the subspace $\mathcal{F}_L$ of Laplace transforms of positive measures. We shall show that within this space, the convolution with a probability measure $g(x)\,dx$ makes a function steeper. (It is shown below that the appropriately shifted $R_\rho * g$ is in $\mathcal{F}_L$.) A key step towards this result, which is also of independent interest, is the following lemma.

LEMMA 7.2. *Let $F = R_{\rho^F} \in \mathcal{F}_L$ satisfy the normalization condition $F(0) = 1$ (i.e., $F = \mathcal{N}F$), and let $G = R_{\rho^G} \in \mathcal{F}_L$ be related to it by*

$$G = \mathcal{N}(F * g) \tag{7.4}$$

*for some probability measure $g(x)\,dx$. Then, for all $\lambda \ge 0$,*

$$\int_0^\lambda \rho^G(du) \le \int_0^\lambda \rho^F(du). \tag{7.5}$$

PROOF. The relation between $F$ and $G$ is such that for some normalizing constant $z \in \mathbb{R}$,

$$G(x) = \int_{-\infty}^{\infty} \left[\int_0^\infty e^{-(x-y+z)u} \rho^F(du)\right] g(y)\,dy = \int_0^\infty e^{-xu} e^{S(u)} \rho^F(du) \tag{7.6}$$

with $S(\cdot)$ defined by

$$e^{S(u)} = \int_{-\infty}^{\infty} e^{(y-z)u} g(y)\,dy. \tag{7.7}$$

Thus

$$\rho^G(du) = e^{S(u)} \rho^F(du), \tag{7.8}$$

and the normalization conditions $F(0) = G(0) = 1$ imply

$$\int_0^\infty e^{S(u)} \rho^F(du) = \int_0^\infty \rho^F(du). \tag{7.9}$$

The function $S(\cdot)$ is convex, which is easily verified by showing that $S'' > 0$, by general properties of integrals of the form (7.7), and satisfies $S(0) = 0$ (since $g(x)\,dx$ is a probability measure). It has, therefore, to be the case that either $\rho^F(du)$ is concentrated at a point (where $S = 0$), or else $S(\cdot)$ is negative on $[0, \bar{u})$ and positive on $(\bar{u}, \infty)$ for some $\bar{u} > 0$. The claimed concentration statement (7.5) is obviously true for all $\lambda \in [0, \bar{u}]$. For $\lambda \ge \bar{u}$, we note that

$$\int_\lambda^\infty e^{S(u)} \rho^F(du) \ge \int_\lambda^\infty \rho^F(du). \tag{7.10}$$

By subtracting (7.9) from (7.10), we find that the claimed (7.5) is valid also for $\lambda > \bar{u}$. $\square$

THEOREM 7.3. *For any $F = R_\rho \in \mathcal{F}_L$ and a probability measure $g(x)\,dx$ on $\mathbb{R}$, the function $\mathcal{N}(F * g)$ is steeper than $F$.*

PROOF. Our goal is to derive the inequality (7.1) for $G = \mathcal{N}(F * g)$ (and $a < b$). By simple approximation arguments, it suffices to do that assuming $\lim_{x\to-\infty} F(x) = \infty$.

We claim that

$$\begin{aligned} &\mathcal{N}(F * g)(x) \le F(x) \qquad \text{for } x \ge 0,\\ &\mathcal{N}(F * g)(x) \ge F(x) \qquad \text{for } x \le 0. \end{aligned} \tag{7.11}$$

We find that the functions $F$ and $G = \mathcal{N}(F * g)$ are related just as in the previous lemma. In order to convert the concentration statement (7.5) into one relating $G(\cdot)$ with $F(\cdot)$, we write, using Fubini's lemma (or integration by parts),

$$\begin{aligned} &\text{for } x > 0 \qquad \mathcal{N}(F * g)(x) = x \int_0^\infty d\lambda e^{-\lambda x}\left[\int_0^\lambda e^{S(u)}\rho(du)\right],\\ &\text{for } x < 0 \qquad \mathcal{N}(F * g)(x) = x \int_0^\infty d\lambda e^{-\lambda x}\left[\int_\lambda^\infty e^{S(u)}\rho(du)\right]\\ &\qquad\qquad\qquad\qquad + \int_0^\infty e^{S(u)}\rho(du), \end{aligned} \tag{7.12}$$

with the corresponding relations holding for $F$ without the factors $e^{S(u)}$. The inequalities (7.11) follow now by inserting here the relations (7.5), (7.10) and (7.9).

We note that if $F$ and $\mathcal{N}(F * g)$ were shifted so as to be equal at a different value of $x$, then the argument above would also go through. Therefore we obtain that $\mathcal{N}(F * g)$ is steeper than $F$. $\square$

The partial order "$G$ is steeper than $F$" is preserved when any of the functions is modified by a uniform shift, and also when each is replaced by a common monotone function of itself, for example, $\{F, G\}$ replaced by $\{1 - e^{-F}, 1 - e^{-G}\}$. Following is a useful property of this partial order (another one is presented in Appendix A.2).

LEMMA 7.4. *Let $F, G \in \mathcal{F}$ be continuous and strictly monotone decreasing functions with $\lim_{x\to-\infty} F(x) = \lim_{x\to-\infty} G(x) = \infty$. If $G$ is steeper than $F$, then, for any $u > 0$,*

$$\int e^{-[G(x-u)-G(x)]}\,de^{-G(x)} \le \int e^{-[F(x-u)-F(x)]}\,de^{-F(x)}. \tag{7.13}$$

*Furthermore, the inequality is strict unless $G$ is a translate of $F$ (and vice versa).*

Proof. The statement is a simple consequence of the following formula, and (7.3):

$$\begin{aligned}&\int e^{-[F(x-u)-F(x)]}\,de^{-F(x)} - \int e^{-[G(x-u)-G(x)]}\,de^{-G(x)}\\ &\quad = \int_0^\infty dz\,[e^{-F(F^{-1}(z)-u)} - e^{-G(G^{-1}(z)-u)}].\end{aligned} \tag{7.14}$$

□

An additional result related to this notion, which may be of independent interest, is presented in Appendix A.2.

**8. Proof of the main result.** We shall now apply the monotonicity arguments for the last leg of the proof of our main result. Theorem 4.2 is clearly implied by the following statement (see Theorem 6.1).

Theorem 8.1. *Let $\mu$ be a measure on the space of configurations $\Omega$, which admits a representation as a random Poisson process, described by a probability measure $\nu(d\rho)$ on $\mathcal{M}$ as in Theorem 6.1, for which both (6.2) and (6.3) hold. Then the support of the Laplace measure $d\rho$ is $\nu$-almost surely a point; that is, the functions $R_\rho$ are almost surely pure exponentials.*

Proof. Let us consider the probability that the first gap exceeds some $u > 0$. For a Poisson process, a simple calculation yields

$$\mathbb{E}_F^{(\text{Poisson})}(x_1 - x_2 \geq u) = \int_{-\infty}^{\infty} e^{-F(x-u)}(-dF(x)). \tag{8.1}$$

Therefore,

$$\mathbb{E}_\mu(x_1 - x_2 \geq u) = \int \mu(d\omega) \int_{-\infty}^{\infty} e^{-F(x-u)}(-dF(x)).$$

Substituting this in (6.2), or in (6.3), one obtains the corresponding expectation for the measure $\mu$. Subtracting the two expressions, we find that

$$0 = \int_{\mathcal{M}} \nu(d\rho)\left[\int_{-\infty}^{\infty} e^{-R_\rho(x-u)}\,dR_\rho(x) - \int_{-\infty}^{\infty} e^{-R_\rho * g(x-u)}\,dR_\rho * g(x)\right]. \tag{8.2}$$

By the analysis in the previous section (Theorem 7.3 and Lemma 7.4), the difference in the square brackets in (8.2) is nonnegative. Thus, this relation implies that

$$\int e^{-R_\rho(x-u)}\,dR_\rho(x) - \int e^{-R_\rho * g(x-u)}\,dR_\rho * g(x) = 0 \qquad \text{for } \nu\text{-almost every } \rho. \tag{8.3}$$

Furthermore, by Lemma 7.4 the equality yields that $\nu$-almost surely $R_\rho * g$ coincides with one of the translates of $R_\rho$. The only functions $(F = R_\rho)$ with this property in $\mathcal{F}_L$ (or for that matter in $\mathcal{F}$; see [3]) are pure exponentials, which correspond to $\rho$ concentrated at a point. □

## APPENDIX

**A.1. Useful statements from the theory of large deviations.** Our goal here is to derive Theorem A.1 which was used in Section 6. Its statement may be read as an expression of the "equivalence of ensembles"—in statistical mechanical terms. The following notation will be used in the theorem.

$$\Lambda(\lambda) \equiv \ln \mathbb{E}[e^{\lambda u_1}], \qquad \Lambda^*(y) \equiv \sup_\lambda (\lambda y - \Lambda(\lambda)).$$

The result we used in Section 6 is:

THEOREM A.1. *Let $u_1, u_2, \dots$ be i.i.d. random variables with a common probability distribution $g$, which has a density and a finite everywhere moment generating function. Then, for any $0 < K < \Lambda'(\infty)$ and $0 < \beta < 1/2$ there is $\varepsilon_{\tau;K,\beta} \underset{\tau\to\infty}{\longrightarrow} 0$ such that for all $q \in [E_g u, K]$ and $|x| \le \tau^\beta$,*

$$\frac{\mathrm{Prob}(\{u_1 + u_2 + \cdots + u_\tau \ge x + q\tau\})}{\mathrm{Prob}(\{u_1 + u_2 + \cdots + u_\tau \ge q\tau\})} = e^{-\eta x}[1 + O(\varepsilon_{\tau;K,\beta})], \tag{A.1}$$

*with $\eta = \eta(q)$ determined by the condition $\eta(q) = \Lambda^{*\prime}(q)$.*

PROOF. We will assume that $E_g u = 0$, since we can replace the random variables $u_i$ by $u_i - E_g u$. We will use the same notation as in the proof of the Bahadur–Rao theorem (see [5]). We denote

$$\eta\left(\frac{y}{\tau}\right) \equiv \Lambda^{*\prime}\left(\frac{y}{\tau}\right),$$

$$\Lambda'\left(\eta\left(\frac{y}{\tau}\right)\right) = \frac{y}{\tau},$$

$$\psi_\tau\left(\eta\left(\frac{y}{\tau}\right)\right) \equiv \eta\left(\frac{y}{\tau}\right)\sqrt{\tau\Lambda''\left(\eta\left(\frac{y}{\tau}\right)\right)},$$

$$Y_i \equiv \frac{u_i - y/\tau}{\sqrt{\Lambda''(\eta(y/\tau))}},$$

$$W_\tau \equiv \frac{Y_1 + \cdots + Y_\tau}{\sqrt{\tau}},$$

and consider a new measure $\widetilde{P}$ defined by its Radon–Nikodym derivative

$$\frac{d\widetilde{P}^{(\eta(y/\tau))}}{dP}(x) = e^{x\eta(y/\tau) - \Lambda(\eta(y/\tau))}.$$

Let also $Q_\tau^{(\eta(y/\tau))}$ denote the distribution function of $W_\tau$ with respect to $\widetilde{P}^{(\eta(y/\tau))}$. It is easy to show then that $Y_i$ are i.i.d. with mean 0 and variance 1 with respect to $\widetilde{P}^{(\eta(y/\tau))}$. Therefore $W_\tau$ has mean 0 and variance 1 with respect to $Q_\tau^{(\eta(y/\tau))}$.

By analogy with the proof of the Bahadur–Rao theorem (see [5]), we can write

$$P(u_1+\cdots+u_\tau \geq y) = e^{-\tau\Lambda^*(y/\tau)} \int_0^\infty e^{-\psi_\tau(\eta(y/\tau))t}\, dQ_\tau^{(\eta(y/\tau))}(t). \tag{A.2}$$

For further consideration, we need to estimate the ratio

$$\begin{aligned} &\frac{P(u_1+u_2+\cdots+u_\tau \geq x+y)}{P(u_1+u_2+\cdots+u_\tau \geq y)} \\ &\qquad = e^{-\tau\Lambda^*((x+y)/\tau)+\tau\Lambda^*(y/\tau)} \frac{\int_0^\infty e^{-\psi_\tau(\eta((x+y)/\tau))t}\, dQ_\tau^{(\eta((x+y)/\tau))}(t)}{\int_0^\infty e^{-\psi_\tau(\eta(y/\tau))t}\, dQ_\tau^{(\eta(y/\tau))}(t)}. \end{aligned} \tag{A.3}$$

By using Taylor's expansion we can estimate the exponent in (A.3):

$$-\Lambda^*\left(\frac{x+y}{\tau}\right) + \Lambda^*\left(\frac{y}{\tau}\right) = -\eta\left(\frac{y}{\tau}\right)\left(\frac{x}{\tau}\right) + O\left(\frac{1}{\tau^{1-2\beta}}\right), \tag{A.4}$$

where, to estimate the remainder term in (A.4), we use the integral form of the remainder in Taylor series and that $\frac{y}{\tau} \leq K$, $|x| \leq \tau^\beta$, $\Lambda^{*\prime\prime} = \frac{1}{\Lambda^{\prime\prime}} < \infty$, convexity of $\Lambda$ and the assumption that the Laplace transform of $g$ is finite.

It remains to show that the prefactor in (A.3) is

$$r(x,y) = \frac{\int_0^\infty e^{-\psi_\tau(\eta((x+y)/\tau))t}\, dQ_\tau^{(\eta((x+y)/\tau))}(t)}{\int_0^\infty e^{-\psi_\tau(\eta(y/\tau))t}\, dQ_\tau^{(\eta(y/\tau))}(t)} = 1 + O(\varepsilon_\tau). \tag{A.5}$$

By the Berry–Esseen theorem (see [7]),

$$\sup_x \left| Q_\tau^{(\eta)}(x) - \int_{-\infty}^x \frac{e^{-t^2/2}}{\sqrt{2\pi}}\, dt \right| \leq \frac{33}{4} \frac{\mathbb{E}u_1^3}{(\operatorname{Var} u_1)^{3/2}} \frac{1}{\sqrt{\tau}} = O\left(\frac{1}{\sqrt{\tau}}\right).$$

Therefore,

$$\begin{aligned} &\int_0^\infty e^{-\psi_\tau(\eta((x+y)/\tau))t}\, dQ_\tau^{(\eta((x+y)/\tau))}(t) \\ &\qquad = \int_0^\infty \frac{e^{-\psi_\tau(\eta((x+y)/\tau))t - t^2/2}}{\sqrt{2\pi}}\, dt \\ &\qquad\quad + O\left(\frac{1}{\sqrt{\tau}}\right)\left(\psi_\tau\left(\eta\left(\frac{x+y}{\tau}\right)\right) + O(1)\right). \end{aligned} \tag{A.6}$$

This formula is especially useful when $\psi_\tau \leq O(1)$ (i.e., when $\eta$ is small) and the first term on the right-hand side of (A.6) is much larger than the second term. In this case we obtain

$$r(x,y)=\frac{\int_0^\infty e^{-\psi_\tau(\eta((x+y)/\tau))t-t^2/2}\,dt+O(1/\sqrt{\tau}\,)}{\int_0^\infty e^{-\psi_\tau(\eta(y/\tau))t-t^2/2}\,dt+O(1/\sqrt{\tau}\,)}. \tag{A.7}$$

If $y$ is such that $O(1)\le\psi_\tau\le O(\tau^{1/2})$, we write the integral as

$$\int_0^\infty e^{-\psi_\tau(\eta(y/\tau))t}\,dQ_\tau^{(\eta(y/\tau))}(t)=\int_0^\infty e^{-\psi_\tau(\eta(y/\tau))t}q_\tau^{(\eta(y/\tau))}(t)\,dt,$$

where $q_\tau$ is the density of $Q_\tau$. By the analog of the Berry–Esseen theorem for densities (see [7]),

$$\sup_x\left|q_\tau(x)-\frac{1}{\sqrt{2\pi}}e^{-x^2/2}\right|=O\left(\frac{1}{\sqrt{\tau}}\right)\qquad\text{as }\tau\to\infty. \tag{A.8}$$

From (A.8) we obtain

$$\begin{aligned} r(x,y)=&\int_0^\infty e^{-\psi_\tau(\eta((x+y)/\tau))t-t^2/2}\,dt+\frac{1}{\psi_\tau(\eta((x+y)/\tau))}O\left(\frac{1}{\sqrt{\tau}}\right)\\ &\times\left\{\int_0^\infty e^{-\psi_\tau(\eta(y/\tau))t-t^2/2}\,dt+\frac{1}{\psi_\tau(\eta(y/\tau))}O\left(\frac{1}{\sqrt{\tau}}\right)\right\}^{-1}. \end{aligned} \tag{A.9}$$

The proof of (A.5) now consists of showing that

$$\begin{aligned} &\int_0^\infty e^{-\psi_\tau(\eta((x+y)/\tau))t-t^2/2}\,dt-\int_0^\infty e^{-\psi_\tau(\eta(y/\tau))t-t^2/2}\,dt\\ &\qquad\le O(\tau^{-\varepsilon})\int_0^\infty e^{-\psi_\tau(\eta(y/\tau))t-t^2/2}\,dt. \end{aligned}$$

□

**A.2. A class of monotone functionals over $\mathcal{F}_L$.** Since the notion introduced in Section 7 may be of independent interest, let us present here a related result, which offers another instructive insight on the contraction properties of convolutions in $\mathcal{F}$.

THEOREM A.2. *Let $F,G\in\mathcal{F}$ with $G$ steeper than $F$. Then, for any positive and continuous function $\Psi$ on $[0,\infty)$ which vanishes at 0 and $\infty$,*

$$\int_{-\infty}^\infty dt\Psi(G(t))\le\int_{-\infty}^\infty dt\Psi(F(t)). \tag{A.10}$$

*Furthermore, if $\Psi$ is strictly positive on $(0,\infty)$, and $G$ and $F$ are both left-continuous, then the inequality is strict unless $G$ is a translate of $F$.*

PROOF. By standard approximation arguments (e.g., using local approximations by polynomials), it suffices to establish (A.10) under the assumption that $\Psi$ is piecewise strictly monotone.

Employing Fubini's lemma, or Lebesgue's "layered cake" formula for the integral,

$$\int_{-\infty}^{\infty} dt \Psi(F(t)) = \int_{0}^{\infty} d\lambda \int_{-\infty}^{\infty} dt I[\Psi(F(t)) \geq \lambda]. \tag{A.11}$$

Under the added assumption on $\Psi$, the set $\{t \in \mathbb{R} : \Psi(F(t)) \geq \lambda\}$ is a union of level intervals of $F$, of the form $\{t \in \mathbb{R} : a_j(\lambda) \leq F(t) \leq b_j(\lambda)\}$ (with $\{[a_j(\lambda), b_j(\lambda)]\}_j$ determined as the level sets of $\{\Psi(\cdot) \geq \lambda\}$).

The integral over $t$ on the right-hand side of (A.11) produces the sum of the lengths of the level-intervals of $F$. When $F$ is replaced by $G$, the corresponding intervals can only get shorter, since $G$ is assumed to be steeper than $F$, and thus (A.10) holds.

In view of the above, the conditions for the *strict* monotonicity sound reasonable. However, since the strict monotonicity is very significant it may be instructive to make the argument explicit. (What follows makes the argument given just above redundant; however, we keep it because of its simplicity.) It is convenient to rearrange the above argument as follows. Using our convention for the inverse function,

$$\int_{-\infty}^{\infty} dt \Psi(F(t)) = \int_{-\infty}^{\infty} dt(F) \Psi(F) = \int_{-\infty}^{\infty} dF^{-1}(a) \Psi(a), \tag{A.12}$$

and thus

$$\begin{aligned} &\int_{-\infty}^{\infty} dt \Psi(F(t)) - \int_{-\infty}^{\infty} dt \Psi(G(t)) \\ &\qquad = \int_{-\infty}^{\infty} [dF^{-1}(a) - dG^{-1}(a)] \Psi(a), \end{aligned} \tag{A.13}$$

where $dF^{-1}(a)$ and $dG^{-1}(a)$ are measures on $\mathbb{R}$. The assumed relation (7.1) implies that the difference $dF^{-1}(a) - dG^{-1}(a)$ is itself a positive measure. The vanishing of its integral against $\Psi$ is therefore possible only if this measure is supported in the set $\Psi^{-1}(0)$, but that set (viewed as the set of values of the functions $F$ and $G$) contains at most the boundary point $a = 0$. It follows that if the inequality (A.10) is saturated, then the two Stieltjes measures are equal in $(0, \infty)$, and thus

$$F^{-1}(a) - G^{-1}(a) = \text{Const}, \tag{A.14}$$

which means that $F$ and $G$ differ by a shift. $\square$

This implies another monotonicity principle, which expresses the fact that convolutions make functions in $\mathcal{F}_L$ steeper.

Corollary A.1. *For any function $F \in \mathcal{F}_L$ and a probability measure $g(x)\,dx$ on $\mathbb{R}$,*

$$\mathbb{E}_{g*F}^{(\text{Poisson})}(x_n - x_{n+1}) \le \mathbb{E}_F^{(\text{Poisson})}(x_n - x_{n+1}) \qquad \textit{for all } n \ge 1, \tag{A.15}$$

*and the inequality is strict unless either both quantities are infinite, or $F(x) = e^{-s(x-z)}$ for some $s > 0$ and $z \in \mathbb{R}$.*

Proof. The mean value of the gap may be computed with the help of the expression

$$x_n - x_{n+1} = \int_{-\infty}^{\infty} \{I[t > x_{n+1}] - I[t > x_n]\}\, dt. \tag{A.16}$$

A simple calculation yields

$$\mathbb{E}_F^{(\text{Poisson})}(x_n - x_{n+1}) = \int_{-\infty}^{\infty} dt \Psi_n(F(t)) \tag{A.17}$$

with $\Psi_n(F) \equiv \frac{F(t)^n}{n!} e^{-F(t)}$. Theorem A.2 applies to such quantities. □

We did not base the proof of Theorem 6.1 on this observation [i.e., use in Section 7 (A.15) instead of (8.1)] since this argument is conclusive only when the above expected value is known to be finite for some $n < \infty$, and we preferred not to limit the proof by such an assumption (and had no need to).

**Acknowledgments.** We thank Pierluigi Contucci for many stimulating discussions in the early part of this project. Anastasia Ruzmaikina would like to express her gratitude to Loren Pitt for his invaluable help and for the financial support from his grant during the course of this work while a Whyburn Research Instructor at the University of Virginia, to Larry Thomas for the reading of an earlier draft of the manuscript and to Almut Burchard, Holly Carley and Etienne DePoortere for useful discussions.

## REFERENCES

[1] Bahadur, R. R. and Rao, R. R. (1960). On deviations of the sample mean. *Ann. Math. Statist.* **31** 1015. MR117775

[2] Bolthausen, E. and Sznitman, A.-S. (1998). On Ruelle's probability cascades and an abstract cavity method. *Comm. Math. Phys.* **197** 247–276. MR1652734

[3] Choquet, G. and Deny, J. (1960). Sur l'équation de convolution $\mu * \sigma = \mu$. *C. R. Acad. Sci. Paris Sér. I Math.* **250** 799–801. MR119041

[4] Daley, D. J. and Vere-Jones, D. (1998). *An Introduction to the Theory of Point Processes*. Springer, Berlin. MR950166

[5] Dembo, A. and Zeitouni, O. (1998). *Large Deviations Techniques and Applications*, 2nd ed. Springer, New York. MR1619036

[6] Derrida, B. (1980). Random-energy model: Limit of a family of disordered models. *Phys. Rev. Lett.* **45** 79–82. MR575260

[7] Feller, W. (1971). *An Introduction to Probability Theory and Its Applications* **2**, 2nd ed. Wiley, New York.

[8] Leadbetter, M. R., Lindgren, G. and Rootzén, H. (1983). *Extremes, and Related Properties of Random Sequences and Processes*. Springer, Berlin. MR691492

[9] Liggett, T. (1979). Random invariant measures for Markov chains, and independent particle systems. *Z. Wahrsch. Verw. Gebiete* **45** 297–854. MR511776

[10] Mezard, M., Parisi, G. and Virasoro, M. A. (1987). *Spin Glass Theory and Beyond*. World Scientific, Singapore. MR1026102

[11] Ruelle, D. (1987). A mathematical reformulation of Derrida's REM and GREM. *Comm. Math. Phys.* **108** 225–239. MR875300

Department of Mathematics
University of Virginia
Charlottesville, Virginia 22903
USA
and
Departments of Statistics
and Mathematics
Purdue University
West Lafayette, Indiana 47905
USA
e-mail: aar@stat.purdue.edu

Departments of Physics
and Mathematics
Princeton University
347 Jadwin Hall
P.O. Box 708
Princeton, New Jersey 08544
USA
e-mail: aizenman@princeton.edu